\documentclass[12pt]{amsart}

\usepackage{amsmath, amsfonts, amsbsy, amsthm, amscd}
\usepackage{amssymb, latexsym, color}
\usepackage{eepic}
\usepackage[dvips]{graphicx}

\newtheorem{Theorem}{Theorem}
\newtheorem{Proposition}[Theorem]{Proposition}
\newtheorem{Lemma}[Theorem]{Lemma}

\newtheorem*{Fact}{Fact}
\newtheorem*{Main}{Main Theorem}

\newtheorem*{Conjecture}{Conjecture}

\newcommand{\C}{{\mathbb{C}}}
\newcommand{\R}{{\mathbb{R}}}

\newcommand{\dps}{\displaystyle}

\newcommand{\area}{\mathrm{Area}}
\newcommand{\ind}{\mathrm{Ind}}
\newcommand{\nul}{\mathrm{Nul}}

\begin{document}

\title{Metrics on a closed surface of genus two which maximize the first 
eigenvalue of the Laplacian}

\author{Shin Nayatani}
\address{Graduate School of Mathematics, Nagoya University, Chikusa-ku, Nagoya 464-8602, Japan}
\email{nayatani@math.nagoya-u.ac.jp}
\thanks{The first author was partly supported by JSPS Grant-in-Aid for Scientific Research (B) 17H02840}

\author{Toshihiro Shoda}
\address{Faculty of Education, Saga University, Honjo-machi, Saga 840-8502, Japan}
\email{tshoda@cc.saga-u.ac.jp}
\thanks{The second author was partly supported by JSPS Grant-in-Aid for Scientific Research (C) 16K05134. }
\thanks{
This work was supported by the Research Institute for Mathematical Sciences, 
a Joint Usage/Research Center located in Kyoto University}


\date{}

\maketitle

\begin{abstract}
In this paper, we settle in the affirmative the Jakobson-Levitin-Nadirashvili-Nigam-Polterovich conjecture, 
stating that a certain singular metric on the Bolza surface, 
with area normalized, should maximize the first eigenvalue of the Laplacian.
\end{abstract} 

\section*{Introduction}

Let $M$ be a closed surface, that is, a compact surface without boundary. 
Throughout this paper, we assume that $M$ is orientable. 
For a Riemannian metric $ds^2$ on $M$, let 
$$
\Lambda(ds^2) := \lambda_1(ds^2)\cdot \area(ds^2), 
$$
where $\lambda_1(ds^2)$ is the first positive eigenvalue of the Laplacian 
and $\area(ds^2)$ is the area of $M$, both with respect to $ds^2$. 
Regarding the upper bound of the quantity $\Lambda(ds^2)$, 
the following results are well-known. 

\begin{Fact}
{\rm (i)\, (Hersch \cite{Hersch})}\,\, For any metric $ds^2$ on the sphere $S^2$, 
$\Lambda(ds^2)\leq 8\pi$ holds.\\ 
{\rm (ii)\, (Yang-Yau \cite{YangYau})}\,\, If $M$ admits a nonconstant meromorphic function 
$(M,ds^2)\rightarrow \overline{\C}=\C\cup \{\infty\}$ of degree $d$, 
then $\Lambda(ds^2)\leq 8\pi\cdot d$ holds. 
In particular, if $\gamma$ is the genus of $M$, then for any metric $ds^2$ on $M$, we have 
\begin{equation}\label{yangyau-ineq}
\Lambda(ds^2)\leq 8\pi\cdot \left[ \frac{\gamma+3}{2} \right]. 
\end{equation}
\end{Fact}

The inequality of (i) is sharp as equality holds
for the standard metric of $S^2$. 
On the other hand, Nadirashvili \cite{Nadirashvili2} found the sharp bound 
$8\pi^2/\sqrt{3}$ of $\Lambda(ds^2)$ for metrics $ds^2$ on the torus $T^2$. 
Thus the inequality \eqref{yangyau-ineq} is not sharp when $\gamma=1$. 



When $\gamma=2$, the inequality \eqref{yangyau-ineq} becomes $\Lambda(ds^2)\leq 16\pi$. 
Jakobson-Levitin-Nadirashvili-Nigam-Polterovich \cite{Nadirashvili1} 
focused their attention on the following metric. 
Let $B$ be the closed Riemann surface of genus two, called the {\em Bolza surface}, 
defined as the smooth completion of the affine complex algebraic curve 
$w^2=z(z^4+1)$. 
Topologically, $B$ is the one-point compactification of the affine curve: 
$$
B = \{ (z,w)\in \C^2 \mid w^2=z(z^4+1) \} \cup \{ (\infty, \infty) \}. 
$$
Let $g_B\colon B\to \overline{\C}$ be the meromorphic function of degree two 
given by $g_B(z,w) = z$. 
If we set $ds^2_B = {g_B}^* ds^2_{S^2}$, where $ds^2_{S^2}$ is the standard metric of 
$S^2=\overline{\C}$, then $ds^2_B$ is a singular Riemannian metric which degenerates 
exactly at the ramification points of $g_B$. 
Since the map $g_B\colon B\rightarrow \overline{\C}$ is a two-sheeted branched 
covering, we have $\area(ds^2_B)=8\pi$. 

\begin{Conjecture}[Jakobson et al.~\cite{Nadirashvili1}]
$\lambda_1(ds^2_B)=2$ should hold. Therefore, $\Lambda(ds^2_B)=16\pi$. 
\end{Conjecture}

For $0< \theta< \pi/2$, let $B_{\theta}$ be the Riemann surface of genus two defined 
as the smooth completion of the affine complex algebraic curve 
$w^2= z (z^4+2\cos 2\theta\cdot z^2+1)$: 
\[
B_{\theta} = \{(z,w)\in \mathbb{C}^2\>| \>
w^2= z (z^4+2\cos 2\theta\cdot z^2+1)\} \cup \{(\infty,\,\infty)\}.
\] 
Note that $B_{\pi/4}=B$. 
Let $ds^2_\theta$ denote the pull-back of the standard metric of $S^2=\overline{\C}$ 
by the meromorphic function $g_\theta\colon B_\theta \ni (z,w) \mapsto z\in \overline{\C}$. 

In this paper, we prove the following theorem, and thereby settle the above conjecture 
in the affirmative. 
\begin{Main}\label{main}
There exists $\theta_1\approx 0.65$ so that for $\theta_1\leq \theta \leq \pi/2 - \theta_1$, 
we have $\lambda_1(ds^2_{B_\theta})=2$ and therefore $\Lambda(ds^2_\theta) = 16\pi$.
\end{Main}

Note that $16\pi$ is a degenerate maximum for $\Lambda$ in the genus two case 
as predicted in \cite{Nadirashvili1}. 
It is also remarked in \cite{Nadirashvili1} that the conjecture implies that
the inequality $\Lambda(ds^2)\leq 16\pi$ is sharp in the class of smooth metrics, 
although the equality may not be attained. 
It is worth mentioning that the Lawson minimal surface of genus two in $S^3$ 
has $\lambda_1 =2$ \cite{ChoeSoret} and $\area \approx 21.91$ \cite{HellerSchmitt}, 
and therefore $\Lambda \approx 43.82<16\pi$.

For recent progress on the existence of $\Lambda$-maximizing metrics 
on a closed surface, see \cite{NadirashviliSire, Petrides}. 

In \S 1, we explain the relation of the above conjecture to the problem of computing 
the Morse index of a minimal surface in Euclidean three-space. 
After that, in \S 2, we prove 
the Main Theorem, assuming two technical lemmas, whose proofs are postponed 
to \S 3 and \S 4.
The paper concludes with two appendices. 

\medskip\noindent
{\em Acknowledgement.}\quad The authors would like to thank Rick Schoen 
\cite{Schoen} for bringing the conjecture of Jakobson et al.~to their attention. 
They would also like to thank the referee for careful reading of the manuscript and 
many invaluable suggestions. 

\section{Index and nullity of a meromorphic function} 

The problem of estimating and computing the Morse (instability) index of a complete 
minimal surface in $\R^3$ (and other flat three-spaces) has been studied by various authors. 
In this section, we explain that the conjecture of Jakobson et al.~is closely 
related to this problem. 

Let $M$ be an orientable complete minimal surface in $\R^3$. 
$M$ is said to be {\em stable} if the second variation of area for any compactly 
supported variation of $M$ is nonnegative, and the plane is the only stable one. 
For non-planar $M$, we define the {\em Morse index} of $M$, $\ind(M)$, as follows: 
For a relatively compact domain $\Omega\subset M$, $\ind(\Omega)$ is defined as 
the maximal dimension of a subspace $V\subset C^\infty_0(\Omega)$ satisfying 
$$
\int_\Omega (|du|^2 + 2Ku^2)\, da<0,\quad \forall u\in V\setminus\{0\}, 
$$
where $K$ and $da$ are the Gaussian curvature and the area element of $M$, 
respectively. 
Note that $\ind(\Omega)$ is necessarily finite. 
We then define 
$$
\ind(M) = \sup_\Omega \ind(\Omega), 
$$
where the supremum is taken over all relatively compact domains $\Omega\subset M$. 
While $\ind(M)$ so defined may become infinity, it was proved by Fischer-Colbrie 
\cite{Fischer-Colbrie} that 
$$
\ind(M)<\infty\quad \Leftrightarrow\quad \int_M (-K)\, da <\infty.
$$ 
Therefore, in studying $\ind(M)$ quantitatively, we may assume that 
$\dps \int_M (-K)\, da <\infty$. In this case, 
$M$ is conformally equivalent to 
a compact Riemann surface $\overline{M}$ with finitely many punctures 
and the Gauss map of $M$, $g\colon M\rightarrow \overline{\C}$, extends to a 
meromorphic function $\overline{g}\colon \overline{M}\rightarrow \overline{\C}$. 
(This is a classical result due to Osserman \cite{Osserman}.) 

In general, for a nonconstant meromorphic function $g\colon M\rightarrow \overline{\C}$
on a compact Riemann surface $M$, we pull back the standard metric of 
$\overline{\C}=S^2$ 
by $g$ and obtain a singular metric $ds^2_g$ (as we did to get $ds^2_B$).
Let $\Delta_g$ denote the Laplacian defined with respect to $ds^2_g$, and 
$\ind(g)$ (resp. $\nul(g)$) 
the number of eigenvalues of $-\Delta_g$ less than $2$ counted with multiplicity 
(resp. the multiplicity of eigenvalue $2$ of $-\Delta_g$). 

\begin{Proposition}[Fischer-Colbrie \cite{Fischer-Colbrie},
Ejiri-Kotani \cite{EjiriKotani}, Montiel-Ros \cite{MontielRos}] 
The Morse index $\ind(M)$ of a complete minimal surface $M$ in $\R^3$ of finite total 
curvature coincides with the index $\ind(\overline{g})$ of the extended Gauss map $\overline{g}$. 
The nullity $\nul(\overline{g})$ equals the dimension of the vector space of 
all bounded Jacobi fields on $M$. 
\end{Proposition}

Since constant functions are necessarily eigenfunctions of the eigenvalue $0$ 
of $-\Delta_g$, we have $\ind(g)\geq 1$. 
The conjecture of Jakobson et al.~asserts that when $g=g_B$, the second least 
eigenvalue 
of $-\Delta_{g_B}$ should equal $2$, and so it is equivalent to asserting 
that $\ind(g_B)=1$. 

\section{Proof of the Main Theorem}

In this section, we prove the Main Theorem, assuming two technical 
Lemmas~{\ref{lemma-extra-ef}} and \ref{lemma-symmetry-extra-ef}. 
The proofs of these lemmas are contained in \S 3 and \S 4.
%
Note that 
the equation of $B_{\theta}$ can be rewritten as 
\[ w^2=z (z-e^{i(\pi/2-\theta)})(z-e^{i(\pi/2+\theta)})(z-e^{-i(\pi/2-\theta)})(z-e^{-i(\pi/2+\theta)}). \] 
Let $g_{\theta}$ and $ds^2_\theta$ be as in the introduction, and $\Delta_\theta$ the Laplacian corresponding to $ds^2_\theta$. 
The meromorphic function $g_{\theta}:B_{\theta}\to \overline{\mathbb{C}}$ gives 
a two-sheeted branched covering which ramifies at the six points 
$(0,\,0),\,(e^{\pm i (\pi/2\pm \theta)},\,0),\,(\infty,\,\infty)$. 
$ds^2_\theta$ is a singular metric which degenerates precisely at the six ramification points 
of $g_\theta$. 
Define three great circular arcs $C_1, C_2, C_3$ on $S^2=\overline{\C}$ by 
\begin{align*}
& C_1 = \{ t \mid t\geq 0 \} \cup \{\infty\},\,\, 
C_2 = \{ e^{i(\pi/2+t)} \mid -\theta \leq t \leq \theta \},\\ 
& C_3 = \{ e^{-i(\pi/2+t)} \mid -\theta \leq t \leq \theta \}. 
\end{align*}
Then $(B_\theta, ds^2_\theta)$ can be represented as the gluing of 
two copies of $(S^2, ds^2_{S^2})$ along $C_1$, $C_2$, $C_3$. 
As $\theta\to 0$, the two arcs $C_2, C_3$ collapse to points, 
and by neglecting the contact at these two points, 
we obtain the metric which is the gluing of two copies of 
$(S^2, ds^2_{S^2})$ along $C_1$. 
The last metric, denoted by $ds^2_0$, is nothing but the pull-back 
of $ds^2_{S^2}$ by the degree two rational function 
$g_0\colon \overline{\C}\ni z\mapsto z^2\in \overline{\C}$. 
Let $\Delta_0$ be the Laplacian defined with respect to $ds^2_0$. 
Then we have the following lemma regarding the eigenvalues 
of $-\Delta_\theta$ and $-\Delta_0$: 

\begin{Lemma}\label{continuity} 
For every positive integer $k$, the $k$-th eigenvalue $\lambda_k(ds^2_\theta)$ of $-\Delta_\theta$ 
is continuous in $\theta$, and as $\theta\to 0$ it converges to the $k$-th eigenvalue 
$\lambda_k(ds^2_0)$ of $-\Delta_0$. 
\end{Lemma}

This lemma may be proved by arguments similar to those in the proof of \cite[Theorem 1]{Nayatani2}. 

In \cite{Nayatani1}, by computing all the eigenvalues of $-\Delta_0$ explicitly, it is shown 
that $\ind(g_0)=3$ and $\nul(g_0)=3$. 
On the other hand, it is known that $\nul(g)\geq 3$ for any nonconstant meromorphic function $g$. 
In fact, the pull-back of three independent eigenfunctions of the eigenvalue $2$ of $-\Delta_{S^2}$, 
the Laplacian with respect to $ds^2_{S^2}$, by $g$ give eigenfunctions of the eigenvalue $2$ 
of $-\Delta_g$. 
From these facts and Lemma \ref{continuity}, it follows that $\ind(g_\theta)=3$ and 
$\nul(g_\theta)=3$ for $\theta$ sufficiently close to $0$. 

We now observe the change of $\nul(g_\theta)$ as $\theta$ increases up to $\pi/4$. 
To do this, we use the work 
of Ejiri-Kotani \cite{EjiriKotani} and Montiel-Ros \cite{MontielRos}. 
If $g$ is a nonconstant meromorphic function such that $\nul(g)>3$, then there exists 
an {\em extra eigenfunction}, that is, an eigenfunction of the eigenvalue $2$ of $-\Delta_g$ 
which is not the pull-back 
of an eigenfunction of the eigenvalue $2$ of $-\Delta_{S^2}$ by $g$. 
As shown in \cite{EjiriKotani, MontielRos}, any extra eigenfunction can be written 
as the support function (that is, the inner product of the position vector field 
and the unit normal vector field) of a complete branched minimal surface 
of finite total curvature whose extended Gauss map is $g$ and whose ends 
are contained in the ramification locus of $g$ and are all planar. 
By using Weierstrass representation, we can express such a minimal surface 
as follows. 
Let $P$ and $B=\sum_{j=1}^l e_j p_j$ 
be the polar and ramification divisors of $g$, respectively, 
where $e_j$ is the multiplicity with which $g$ takes its value at $p_j$. 
Set $D=B-2P$. 
Suppose that there exists a non-zero $\omega\in H^0 (M,\, K_M\otimes D)$ satisfying 
\begin{equation}\label{NS-residue}
{\rm Res}_{p_j}\omega=0,\quad 1\leq \forall j\leq l, 
\end{equation}
and 
\begin{equation}\label{NS-period}
\Re\int_{\ell} {}^t (1-g^2,\,i(1+g^2),\,2g) \omega = \textbf{o},\quad 
\forall \ell\in H_1(M,\,\mathbb{Z}),
\end{equation}
where $K_M$ is the canonical divisor of $M$. 
Then for any such $\omega$, 
\[
X_{\omega} (p) =\Re\int^p_{p_0} {}^t (1-g^2,\,i(1+g^2),\,2g)\, \omega 
\]
gives a minimal surface with the above properties. 
%

We now apply the general result as above to $(B_\theta, g_\theta)$. 
We can determine the values of $\theta$ for which there exists a non-zero 
$\omega\in H^0 (B_\theta,\, K_{B_\theta}\otimes D)$ satisfying 
\eqref{NS-residue} and \eqref{NS-period}. 
In fact, we have
\begin{Lemma}\label{lemma-extra-ef}
Set 
\begin{align*}
&A=\int_0^{\infty} \frac{d t}{\sqrt{t(t^4 +2 \cos 2\theta\cdot t^2+1)}}, \quad 
B=\int_0^{\infty} \frac{d t}{\sqrt{t(t^4 -2 \cos 2\theta\cdot  t^2+1)}}, \\
&C= \int_0^{\infty} \frac{t^3\,d t}{\sqrt{t(t^4 +2 \cos 2\theta\cdot  t^2+1)}^3},\quad 
D= \int_0^{\infty} \frac{t^3\,d t}{\sqrt{t(t^4 -2 \cos 2\theta\cdot  t^2+1)}^3}.
\end{align*}
Let $\theta_1$\,{\rm (}$\approx 0.65${\rm )} be the unique solution of 
\[
A (B^2+16 D^2\sin^2 2\theta) +8 (AD+BC)(B \cos 2\theta - 4 D \sin^2 2\theta)=0,
\] 
and set $\theta_2=\pi/2-\theta_1$\,{\rm (}$\approx 0.91${\rm )}. 
Then there exists a non-zero 
$\omega\in H^0 (B_\theta,\, K_{B_\theta}\otimes D)$ satisfying 
\eqref{NS-residue} and \eqref{NS-period} if and only if $\theta=\theta_1,\,\theta_2$. 
If $\theta=\theta_1$, then any such $\omega$ is given by a real linear combination of 
\begin{align*}
\omega_1 &:=  -\frac{AD+3 BC}{4(AD+BC)} \frac{dz}{w}
-\frac{AD+3 BC}{4(AD+BC)} \frac{dz}{w^3} +\frac{z}{w^3}dz \\
&+\frac{AB+(AD-BC)\cos 2\theta}{2(AD+BC)} \frac{z^2}{w^3}dz
+\frac{AB+2(AD+BC)\cos 2\theta}{2(AD+BC)} \frac{z^3}{w^3}dz \\
&+\frac{3AD+ BC}{4(AD+BC)} \frac{z^4}{w^3}dz, \\
\omega_2 &:= i \, \bigg( -\frac{AD+3 BC}{4(AD+BC)} \frac{dz}{w}
+\frac{AD+3 BC}{4(AD+BC)} \frac{dz}{w^3} +\frac{z}{w^3}dz \\
&-\frac{AB+(AD-BC)\cos 2\theta}{2(AD+BC)} \frac{z^2}{w^3}dz
+\frac{AB+2(AD+BC)\cos 2\theta}{2(AD+BC)} \frac{z^3}{w^3}dz \\
&-\frac{3AD+ BC}{4(AD+BC)} \frac{z^4}{w^3}dz \bigg). 
\end{align*}
{\rm (}We can obtain a similar assertion for $\theta=\theta_2$.{\rm )}
\end{Lemma}

The lemma implies that there are two independent extra eigenfunctions 
when $\theta = \theta_1,\,\theta_2$. 
Thus we obtain 
\begin{Proposition}
\begin{equation}
{\rm Nul}\, (g_\theta)=
\begin{cases}
5, & \theta=\theta_1,\,\theta_2, \\
3, & \theta\neq \theta_1,\,\theta_2.
\end{cases}
\end{equation}
\end{Proposition}

To see how ${\rm Ind}\, (g_{\theta})$ changes as $\theta$ increases and 
passes $\theta_1$, we use symmetries of $B_{\theta}$. 
Let $j:B_{\theta}\to B_{\theta}$ be the hyperelliptic involution given by $j(z,w)=(z,-w)$, 
and $s_1,\,s_2,\,s_3:B_{\theta}\to B_{\theta}$ the anti-holomorphic involutions given by 
$s_1 (z,w)=(\overline{z},\overline{w}),\, s_2 (z,w)=(-\overline{z},i\,\overline{w}),\, 
s_3 (z,w)=(1/\overline{z},\overline{w}/\overline{z}^3)$. 
We have 
\[
s_1\circ s_2 = j\circ s_2\circ s_1,\quad s_2\circ s_3=s_3\circ s_2,\quad s_3\circ s_1=s_1\circ s_3.
\]
Thus the three involutions $j,\,s_1,\,s_3$ of $B_{\theta}$ commute with one another, and 
the group of symmetries, $H$, generated by them is an abelian group 
of order eight. 
A fundamental domain for the action of $H$ on $B_{\theta}$ is given by 
the intersection of the upper half plane and the unit disk, 
denoted by $\Omega$. (See Figure \ref{NS-fundamental-domain}.) 

\begin{figure}[htbp]
\begin{picture}(340,120)

\put(200,95){\circle*{3}} \put(205,95){$e^{i(\pi/2-\theta)}$}
\put(140,95){\circle*{3}} \put(100,90){$e^{i(\pi/2+\theta)}$}
\put(200,15){\circle*{3}} \put(205,10){$e^{-i(\pi/2-\theta)}$}
\put(140,15){\circle*{3}} \put(92,10){$e^{-i(\pi/2+\theta)}$}
\put(170,55){\circle*{3}} \put(172,43){$O$}

\thicklines

\textcolor{blue}{\qbezier(140,95)(170,115)(200,95)}
\!\!\textcolor{blue}{\qbezier(140,15)(170,-5)(200,15)}

\!\!\textcolor{red}{\qbezier(140,95)(120,79)(120,55)}
\!\!\textcolor{red}{\qbezier(140,15)(120,31)(120,55)}
\!\!\textcolor{red}{\qbezier(200,95)(220,79)(220,55)}
\!\!\textcolor{red}{\qbezier(200,15)(220,31)(220,55)}

\put(170,55){\textcolor{red}{\vector(1,0){70}}} 
\put(170,55){\textcolor{blue}{\line(-1,0){70}}} 
\put(170,-8){\vector(0,1){130}}

\end{picture}
\caption{Fixed point sets of $s_1,\,j\circ s_1,\, s_3,\,j\circ s_3$ 
}\label{NS-fixed-point-sets}
\end{figure}
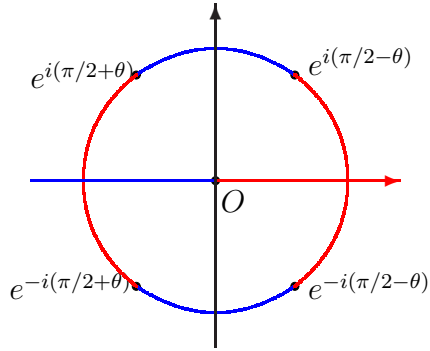

Recall that $B_\theta$ is the gluing of two copies of $\overline{\mathbb{C}}$. 
The fixed point sets of the anti-holomorphic involutions 
$s_1,\,j\circ s_1,\, s_3,\,j\circ s_3$ are as follows. 
(See Figure \ref{NS-fixed-point-sets}.) 
\begin{itemize}
\item The fixed point set of $s_1$ is the red half-line on the real axis, 
\item The fixed point set of $j\circ s_1$ is the blue half-line on the real axis, 
\item The fixed point set of $s_3$ is the union of the red arcs on the unit circle, 
\item The fixed point set of $j\circ s_3$ is the union of the blue arcs on the unit circle. 
\end{itemize}
For example, $s_1 (z,\,w)=(z,\,w)$ if and only if $z$ ($=:x$), $w$ ($=:y$) are real. 
Since 
$$
y^2 = x(x^4+2\cos 2\theta\cdot x^2+1) = x\{(x^2+\cos 2\theta)^2+\sin^2 2\theta\}\geq 0
$$
and $(x^2+\cos 2\theta)^2+\sin^2 2\theta> 0$, one must have $x\geq 0$. 

\begin{figure}[htbp]
\begin{picture}(340,78)

\put(177,30){$\Omega$}

\put(200,50){\circle*{3}} \put(205,50){$e^{i (\pi/2-\theta)}$}
\put(140,50){\circle*{3}} \put(100,45){$e^{i (\pi/2+\theta)}$}
\put(170,10){\circle*{3}} \put(173,-2){$O$}

\thicklines

\qbezier(140,50)(170,70)(200,50)

\qbezier(140,50)(120,34)(120,10)
\qbezier(200,50)(220,34)(220,10)

\put(170,10){\line(1,0){50}} 
\put(170,10){\line(-1,0){50}} 
\put(120,10){\line(-1,0){20}} \put(220,10){\vector(1,0){20}}
\put(170,-10){\vector(0,1){90}}

\end{picture}
\caption{The fundamental domain $\Omega$ for $H$}\label{NS-fundamental-domain}
\end{figure}
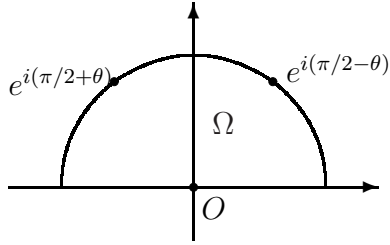

Since $H$ is abelian and preserves $ds^2_\theta$, each eigenspace of 
$-\Delta_\theta$ is invariant under the action of $H$ and 
spanned by simultaneous eigenvectors for all $s\in H$. 
Let $u_i$, $i=1,\,2$, be the support functions of the branched minimal 
immersions $X_{\omega_i}$, in whose definition we choose $p_0=(1, 
\sqrt{2+2\cos 2\theta})$ as the base point. 
They are extra eigenfunctions for $\theta=\theta_1$. 
The following lemma shows how $H$ acts on $u_1,\,u_2$. 
\begin{Lemma}\label{lemma-symmetry-extra-ef}
\begin{align*}
s_1^{*} u_1&=u_1,\quad s_3^{*} u_1=u_1,\quad 
j^{*} u_1=-u_1 + \langle c_1, N \rangle, \\
s_1^{*} u_2&=-u_2,\quad s_3^{*} u_2=u_2,\quad 
j^{*} u_2=-u_2 + \langle c_2, N \rangle, 
\end{align*}
where $c_i\in \R^3$, $i=1,2$, and $N$ is the unit normal vector field 
of $X_{\omega_i}$. 
\end{Lemma}
In order to get extra eigenfunctions which behave properly with respect to 
the actions of $j\circ s_1$ and $j\circ s_3$, we set
\begin{align*}
&v_1 = u_1-(j\circ s_1)^{*} u_1-(j\circ s_3)^{*} u_1+
(j\circ s_1)^{*}\circ (j\circ s_3)^{*} u_1, \\
&v_2 = u_2+(j\circ s_1)^{*} u_2-(j\circ s_3)^{*} u_2-
(j\circ s_1)^{*}\circ (j\circ s_3)^{*} u_2. 
\end{align*}
By Lemma~{\ref{lemma-symmetry-extra-ef}}, we have 
\begin{align*}
&s_1^{*} v_1=v_1,\quad (j\circ s_1)^{*} v_1=-v_1,\quad 
s_3^{*} v_1=v_1, \quad  (j\circ s_3)^{*} v_1=-v_1, \\
& s_1^{*} v_2=-v_2,\quad (j\circ s_1)^{*} v_2=v_2,\quad 
s_3^{*} v_2=v_2,\quad (j\circ s_3)^{*} v_2=-v_2. 
\end{align*}

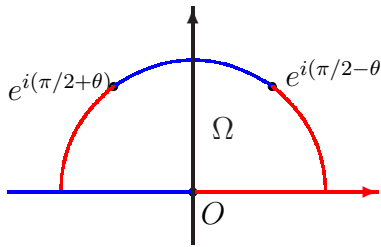
\begin{figure}[htbp]
\begin{picture}(340,75)

\put(177,30){$\Omega$}

\put(200,50){\circle*{3}} \put(205,50){$e^{i (\pi/2-\theta)}$}
\put(140,50){\circle*{3}} \put(100,45){$e^{i (\pi/2+\theta)}$}
\put(170,10){\circle*{3}} \put(173,-2){$O$}

\thicklines

\textcolor{blue}{\qbezier(140,50)(170,70)(200,50)}

\!\!\textcolor{red}{\qbezier(140,50)(120,34)(120,10)}
\!\!\textcolor{red}{\qbezier(200,50)(220,34)(220,10)}

\put(170,10){\textcolor{red}{\line(1,0){50}}} 
\put(170,10){\textcolor{blue}{\line(-1,0){70}}} 
\put(220,10){\textcolor{red}{\vector(1,0){20}}}
\put(170,-10){\vector(0,1){90}}

\end{picture}
\caption{Fixed point sets in $\partial \Omega$}\label{NS-fixed-point-sets2}
\end{figure}

Henceforth, we regard $v_1$ and $v_2$ 
as functions on $\Omega$. (See Figure \ref{NS-fixed-point-sets2}.)
Then the preceding observations mean that $v_1$ 
satisfies the Dirichlet (resp. Neumann) 
condition on the blue (resp. red) segments in the unit circle 
and on the blue (resp. red) segment in the real axis. 
As $\theta$ increases, 
the blue (resp. red) segment in the unit circle becomes longer (resp. shorter). 
Hence, by the variational characterization of eigenvalues, the eigenvalues of the Laplacian 
in $\Omega$ under the boundary conditions as above monotonically increase. 
Similarly, $v_2$ satisfies the Dirichlet (resp. Neumann) condition on the blue (resp. red) 
segment in the unit circle and on the red (resp. blue) segment in the real axis, 
and therefore the eigenvalues of the Laplacian in $\Omega$ under the boundary conditions of $v_2$ 
also monotonically increase. 

The two assertions we just made mean that there exist two independent eigenfunctions 
of $-\Delta_\theta$ with the same type of symmetry as $v_1$ 
and $v_2$, respectively, 
such that the corresponding eigenvalues increase monotonically and continuously. 
On the other hand, for $0< \theta< \theta_2$, extra eigenfunctions with the other types 
of symmetry do not occur. 
Hence, the number of the eigenvalues of $-\Delta_\theta$ less than $2$, 
whose eigenfunctions have the other types of symmetry, 
remains unchanged throughout $(0, \theta_2)$. 
(Here we use the continuity of eigenvalues in $\theta$ again.) 

We may now conclude that as $\theta$ increases and passes 
$\theta_1$, two eigenvalues of 
$-\Delta_\theta$ will monotonically increase and pass $2$ upward, 
and thus the number of eigenvalues less than $2$ decreases by two. 
One can also verify that if $\theta$ increases further and passes 
$\theta_2$, then two eigenvalues of 
$-\Delta_\theta$ will decrease and pass $2$ downward, and the number of eigenvalues 
less than $2$ increases by two. 
To summarize, we have proved the following 

\begin{Theorem}\label{theorem-index}
\begin{equation}\label{index}
\ind(g_\theta) = \left\{\begin{array}{cl}
3, & 0 < \theta < \theta_1,\\ 
1, & \theta_1\leq \theta \leq \theta_2,\\ 
3, & \theta_2 < \theta < \pi/2.
\end{array}\right.
\end{equation}
\end{Theorem}

This theorem implies the Main Theorem. 

\section{Proof of Lemma \ref{lemma-extra-ef}}

This section is devoted to the proof of Lemma \ref{lemma-extra-ef}. 

Recall that $K_{B_\theta}$ is the canonical divisor of $B_\theta$ and $D=B-2P$, 
where $P$ and $B=\sum_{j=1}^l e_j p_j$ are the polar and ramification divisors 
of $g_\theta$, respectively.
Let $\widehat{H}(g_\theta)$ denote the set of all 
$\omega\in H^0 (B_\theta,\, K_{B_\theta}\otimes D)$ satisfying 
\begin{equation}\label{residue}
{\rm Res}_{p_j}\omega=0,\quad 1\leq \forall j\leq l, 
\end{equation}
and $H(g_\theta)$ the set of all $\omega\in \widehat{H}(g_\theta)$ satisfying 
\begin{equation}\label{period}
\Re\int_{\ell} {}^t\! \left( 1-{g_\theta}^2,\,i(1+{g_\theta}^2),\,2g_\theta \right)\, 
\omega = \textbf{o},\quad \forall \ell\in H_1(B_{\theta},\,\mathbb{Z}). 
\end{equation}
Note that $\widehat{H}(g_\theta)$ is a complex vector space. 
We should determine the values of $\theta$ for which $H(g_\theta)\neq \{0\}$. 

We first find a basis for $\widehat{H}(g_\theta)$. 
The polar and ramification divisors of $g_{\theta}$ are given by 
$$
P=2(\infty,\,\infty),\quad 
B = 2(0,\,0)+2(e^{\pm i (\pi/2\pm\theta)},\,0)+2(\infty,\,\infty), 
$$
and therefore  
$$
D = 2(0,\,0)+2(e^{\pm i (\pi/2\pm\theta)},\,0)-2(\infty,\,\infty). 
$$
By the Riemann-Roch theorem, $H^0 (B_{\theta}, K_{B_{\theta}}\otimes D)$ 
has dimension nine, and 
$$
\left\{ \frac{dz}{w},\, \frac{dz}{w^2},\,\frac{z}{w^2} dz,\, \frac{z^2}{w^2} dz,\, 
\frac{dz}{w^3},\, \frac{z}{w^3} dz,\, \frac{z^2}{w^3} dz,\, \frac{z^3}{w^3} dz,\, 
\frac{z^4}{w^3} dz \right\} 
$$
is a basis for it. 
It is easy to verify that 
$$
\left\{ \frac{dz}{w},\, \frac{dz}{w^3},\, \frac{z}{w^3} dz,\, \frac{z^2}{w^3} dz,\, 
\frac{z^3}{w^3} dz,\, \frac{z^4}{w^3} dz \right\} 
$$ 
is a basis for $\widehat{H}(g_\theta)$. 
Therefore, $\omega\in \widehat{H}(g_\theta)$ has the form 
\begin{equation}\label{linearcombination}
\omega=\alpha_1 \frac{dz}{w} + \alpha_2 \frac{dz}{w^3}+\alpha_3 \frac{z}{w^3} dz 
+ \alpha_4\frac{z^2}{w^3} dz+ \alpha_5 \frac{z^3}{w^3} dz+\alpha_6 \frac{z^4}{w^3} dz, 
\end{equation}
where $\alpha_1,\,\ldots,\,\alpha_6$ are complex numbers. 

We now consider the period condition \eqref{period}. 
First, we express the above basis elements of $\widehat{H}(g_\theta)$ 
as linear combinations of the abelian differentials of the second kind 
$dz/w, zdz/w, z^3 dz/w^3, z^4 dz/w^3$ up to exact forms. 
It is easy to show that
\begin{align}
d (z^p w^q) &= \frac{1}{2} z^{p-1} w^{q-2} \{ (2p+5q)w^2
-4q \cos 2\theta\cdot  z^3-4qz\} dz \label{Nayatani-Shoda-exact1}\\
&= \frac{1}{2} z^p w^{q-2} \{ (2p+5q)z^4+2\cos 2\theta\cdot (2p+3q) z^2 
+2p+q \} dz.
\label{Nayatani-Shoda-exact2}
\end{align}
For two meromorphic one-forms $\eta_1$, $\eta_2$ on $B_{\theta}$, we write 
$\eta_1\sim \eta_2$ if there exists a meromorphic function $f$ on $B_{\theta}$ 
such that $\eta_1=\eta_2+d f$.
By using \eqref{Nayatani-Shoda-exact1}, \eqref{Nayatani-Shoda-exact2} we deduce 
the following relations: 
\begin{align}
& \frac{z}{w^3} dz \sim \frac{3}{4} \frac{dz}{w} - \cos 2\theta\cdot  \frac{z^3}{w^3} dz, 
\label{Nayatani-Shoda-key1}\\
& \frac{z^2}{w^3} dz \sim \frac{1}{4} \frac{z}{w} dz - \cos 2\theta\cdot  \frac{z^4}{w^3} dz, 
\label{Nayatani-Shoda-key2} \\
& \frac{dz}{w^3} \sim -\frac{3}{2} \cos 2\theta\cdot  \frac{z}{w} dz
+(-5+6 \cos^2 2\theta) \frac{z^4}{w^3} dz, \label{Nayatani-Shoda-key3} \\
& \frac{z^5}{w^3} dz \sim \frac{1}{4} \frac{dz}{w}-\cos 2\theta\cdot  \frac{z^3}{w^3} dz, \label{Nayatani-Shoda-key4} \\
& \frac{z^6}{w^3} dz 
\sim \frac{3}{4} \frac{z}{w} dz-\cos 2\theta\cdot  \frac{z^4}{w^3} dz, 
\label{Nayatani-Shoda-key5} \\
& \frac{z^2}{w}dz 
\sim -\cos 2\theta\cdot  \frac{dz}{w}-4\sin^2 2\theta\cdot  \frac{z^3}{w^3}dz. 
\label{Nayatani-Shoda-key6}
\end{align}
In fact, \eqref{Nayatani-Shoda-key1} and \eqref{Nayatani-Shoda-key2} follow immediately
from \eqref{Nayatani-Shoda-exact1} with choices $(p,q)=(1,-1)$ and $(p,q)=(2,-1)$, 
respectively. 
\eqref{Nayatani-Shoda-key3} follows by using \eqref{Nayatani-Shoda-exact2} with 
$(p,q)=(0,-1)$ and then applying \eqref{Nayatani-Shoda-key2}.
\eqref{Nayatani-Shoda-key4} and \eqref{Nayatani-Shoda-key5} follow by substituting 
$z^5 = w^2-2\cos2\theta\cdot z^3-z$ and then applying \eqref{Nayatani-Shoda-key1} and 
\eqref{Nayatani-Shoda-key2}, respectively. 
Finally, \eqref{Nayatani-Shoda-key6} follows by using \eqref{Nayatani-Shoda-exact1} 
with $(p,q)=(3,-1)$ and then applying \eqref{Nayatani-Shoda-key4}. 

For 
\[
\omega=\alpha_1 \frac{dz}{w} + \alpha_2 \frac{dz}{w^3}+\alpha_3 \frac{z}{w^3} dz 
+ \alpha_4\frac{z^2}{w^3} dz+ \alpha_5 \frac{z^3}{w^3} dz+\alpha_6 \frac{z^4}{w^3} dz 
\]
as in \eqref{linearcombination}, we find by using the above relations
\begin{align}\label{Nayatani-Shoda-period-cond2} 
\omega \underbrace{\sim}_{\eqref{Nayatani-Shoda-key1}, \,\eqref{Nayatani-Shoda-key2}, \, \eqref{Nayatani-Shoda-key3}}& 
\left(\alpha_1 + \frac{3}{4} \alpha_3 \right) \frac{dz}{w} 
+\left(-\frac{3}{2} \cos 2\theta\cdot  \alpha_2 +\frac{\alpha_4}{4} \right) \frac{z}{w} dz 
\\ 
& +  (-\cos 2\theta\cdot  \alpha_3 +\alpha_5 )\frac{z^3}{w^3} dz \nonumber \\
& +((-5+6 \cos^2 2\theta)\alpha_2 -\cos 2\theta\cdot  \alpha_4 + \alpha_6) \frac{z^4}{w^3} dz, 
\nonumber \\
\label{Nayatani-Shoda-period-cond3} z \omega 
\underbrace{\sim}_{\eqref{Nayatani-Shoda-key1}, \,\eqref{Nayatani-Shoda-key2}, \,
\eqref{Nayatani-Shoda-key4}}& 
\left(\frac{3}{4}\alpha_2 + \frac{\alpha_6}{4} \right) \frac{dz}{w} 
+ \left(\alpha_1 + \frac{\alpha_3}{4}\right) \frac{z}{w} dz\\ 
& + (-\cos 2\theta\cdot  \alpha_2+\alpha_4-\cos 2\theta\cdot  \alpha_6)\frac{z^3}{w^3} dz 
\nonumber \\ 
& + (-\cos 2\theta\cdot  \alpha_3+\alpha_5) \frac{z^4}{w^3} dz, \nonumber
\end{align}
and 
\begin{align}\label{Nayatani-Shoda-period-cond4} 
z^2 \omega 
\underbrace{\sim}_{\eqref{Nayatani-Shoda-key2}, \,\eqref{Nayatani-Shoda-key4}, \,
\eqref{Nayatani-Shoda-key5},\,\eqref{Nayatani-Shoda-key6}}& 
\left(-\cos 2\theta\cdot  \alpha_1 +\frac{\alpha_5}{4} \right) \frac{dz}{w} 
+ \left(\frac{\alpha_2}{4}+\frac{3}{4} \alpha_6 \right) \frac{z}{w} dz
\end{align}
\begin{align}
& + (-4\sin^2 2\theta\cdot  \alpha_1 + \alpha_3-\cos 2\theta\cdot  \alpha_5) 
\frac{z^3}{w^3}dz \nonumber \\
& + (-\cos 2\theta\cdot  \alpha_2+\alpha_4-\cos 2\theta\cdot  \alpha_6)\frac{z^4}{w^3} dz. \nonumber 
\end{align}

Let 
$\varphi\colon B_\theta\to B_\theta$ be 
the automorphism given by $\varphi (z,\,w)=(-z,\,i w)$. 
Note that $\varphi^2=j$, the hyperelliptic involution of $B_\theta$. 
Define paths $C_4$, $C_5$ on $B_\theta$ by 
\begin{align*}
C_4 &= \{(z,\,w)=(t,\,\sqrt{t(t^4 +2 \cos 2\theta\cdot  t^2+1)}\,)\>|\> 0\leq t\leq \infty\}, \\
C_5 &= \{(z,\,w)=(i t,\,e^{i \pi/4}\sqrt{t(t^4 -2 \cos 2\theta\cdot  t^2+1)}\,)\>|\> 0\leq t\leq \infty\}. 
\end{align*}
Then the four closed paths 
\begin{align*}
C_4\cup (-j(C_4)),\> \varphi (C_4\cup (-j(C_4))), \> 
C_5\cup (-j(C_5)),\> \varphi (C_5\cup (-j(C_5))) 
\end{align*}
form a homology basis, as verified by integrating the holomorphic 
differentials $dz/w,\>z dz/w$ over them. 

Straightforward calculations yield 
\begin{align*}
&\int_{C_4\cup \{-j(C_4)\}} \frac{dz}{w} 
= 2 \int_0^{\infty} \frac{d t}{\sqrt{t(t^4 +2 \cos 2\theta\cdot  t^2+1)}}= 2 A,\\ 
&\int_{C_4\cup \{-j(C_4)\}} \frac{z}{w} dz 
= 2 \int_0^{\infty} \frac{t\, d t}{\sqrt{t(t^4 +2 \cos 2\theta\cdot  t^2+1)}} 
\underbrace{=}_{s=1/t} 2A, \\
&\int_{\varphi ( C_4\cup \{-j(C_4)\})} \frac{dz}{w} 
= 2\, i A,\quad  
\int_{\varphi (C_4\cup \{-j(C_4)\})} \frac{z}{w} dz 
= -2\, i A, \\
&\int_{C_5\cup \{-j(C_5)\}} \frac{dz}{w} 
= 2\, e^{\frac{\pi}{4} i}B,\quad
\int_{C_5\cup \{-j(C_5)\}} \frac{z}{w} dz 
= -2\, e^{-\frac{\pi}{4} i}B, \\
&\int_{\varphi ( C_5\cup \{-j(C_5)\})} \frac{dz}{w} 
= -2\, e^{-\frac{\pi}{4} i} B,\quad 
\int_{\varphi (C_5\cup \{-j(C_5)\})} \frac{z}{w} dz 
= 2\, e^{\frac{\pi}{4} i} B,\\ 
&\int_{C_4\cup \{-j(C_4)\}} \frac{z^3}{w^3} dz
= 2 \int_0^{\infty} \frac{t^3\,d t}{\sqrt{t(t^4 +2 \cos 2\theta\cdot  t^2+1)}^3}= 2\, C,\\ 
&\int_{C_4\cup \{-j(C_4)\}} \frac{z^4}{w^3} dz
= 2 \int_0^{\infty} \frac{t^4\,d t}{\sqrt{t(t^4 +2 \cos 2\theta\cdot  t^2+1)}^3}
\underbrace{=}_{s=1/t} 2\,C, \\
&\int_{\varphi ( C_4\cup \{-j(C_4)\})} \frac{z^3}{w^3} dz
= 2\, i \, C,\quad 
\int_{\varphi (C_4\cup \{-j(C_4)\})} \frac{z^4}{w^3} dz 
= -2\, i \, C, \\
&\int_{C_5\cup \{-j(C_5)\}} \frac{z^3}{w^3} dz
= -2\, e^{\frac{\pi}{4} i} D,\quad 
\int_{C_5\cup \{-j(C_5)\}} \frac{z^4}{w^3} dz
= 2\, e^{-\frac{\pi}{4} i}D, \\
&\int_{\varphi ( C_5\cup \{-j(C_5)\})} \frac{z^3}{w^3} dz
= 2\, e^{-\frac{\pi}{4} i} D,\quad
\int_{\varphi (C_5\cup \{-j(C_5)\})} \frac{z^4}{w^3} dz
= -2\, e^{\frac{\pi}{4} i} D. 
\end{align*}

Note that the period condition \eqref{period} can be rewritten as 
\begin{equation}\label{period-new}
\int_{\ell}  \omega = \overline{\int_{\ell} g_\theta^2\, \omega},\quad 
\Re \int_{\ell} g_\theta\, \omega = \textbf{o},\quad 
\forall \ell\in H_1(B_{\theta},\,\mathbb{Z}).
\end{equation}
By using \eqref{Nayatani-Shoda-period-cond2}\! --\! \eqref{Nayatani-Shoda-period-cond4} and 
the calculation we have just made, one can express the former relation of \eqref{period-new} 
for the above homology basis as 
%
%
\begin{align}
\label{Nayatani-Shoda-eq1} &\left(\alpha_1 + \frac{3}{4} \alpha_3 \right) A
+\left(-\frac{3}{2} \cos 2\theta\cdot  \alpha_2 +\frac{\alpha_4}{4} \right) A\\ 
&\quad +  (-\cos 2\theta\cdot  \alpha_3 +\alpha_5 ) C   
+((-5+6 \cos^2 2\theta)\alpha_2 -\cos 2\theta\cdot  \alpha_4 + \alpha_6) C \nonumber \\
&= 
\overline{\left(-\cos 2\theta\cdot  \alpha_1 +\frac{\alpha_5}{4} \right) A 
+ \left(\frac{\alpha_2}{4}+\frac{3}{4} \alpha_6 \right) A} \nonumber \\
&\quad \overline{+ (-4\sin^2 2\theta\cdot  \alpha_1 + \alpha_3-\cos 2\theta\cdot  \alpha_5) C 
+ (-\cos 2\theta\cdot  \alpha_2+\alpha_4-\cos 2\theta\cdot  \alpha_6) C}, \nonumber \\
\label{Nayatani-Shoda-eq3} &\left(\alpha_1 + \frac{3}{4} \alpha_3 \right) i A
+\left(-\frac{3}{2} \cos 2\theta\cdot  \alpha_2 +\frac{\alpha_4}{4} \right) (-i A) \\ 
&\quad +  (-\cos 2\theta\cdot  \alpha_3 +\alpha_5 ) i C   
+((-5+6 \cos^2 2\theta)\alpha_2 -\cos 2\theta\cdot  \alpha_4 + \alpha_6) (-i C) \nonumber \\
&= 
\overline{\left(-\cos 2\theta\cdot  \alpha_1 +\frac{\alpha_5}{4} \right) i A 
+ \left(\frac{\alpha_2}{4}+\frac{3}{4} \alpha_6 \right) (-i A) } \nonumber \\
&\quad \overline{+ (-4\sin^2 2\theta\cdot  \alpha_1 + \alpha_3-\cos 2\theta\cdot  \alpha_5) i\, C 
+ (-\cos 2\theta\cdot  \alpha_2+\alpha_4-\cos 2\theta\cdot  \alpha_6) (-i C)}, \nonumber \\
\label{Nayatani-Shoda-eq5} &\left(\alpha_1 + \frac{3}{4} \alpha_3 \right) (1+i) B
+\left(-\frac{3}{2} \cos 2\theta\cdot  \alpha_2 +\frac{\alpha_4}{4} \right) (-1+i) B \\ 
&\qquad +  (-\cos 2\theta\cdot  \alpha_3 +\alpha_5 ) (-1-i) D \nonumber \\ 
&\qquad +((-5+6 \cos^2 2\theta)\alpha_2 -\cos 2\theta\cdot  \alpha_4 + \alpha_6) (1-i) D 
\nonumber \\
&= 
\overline{\left(-\cos 2\theta\cdot  \alpha_1 +\frac{\alpha_5}{4} \right) (1+i) B 
+ \left(\frac{\alpha_2}{4}+\frac{3}{4} \alpha_6 \right) (-1+i) B} \nonumber \\
&\qquad \overline{+ (-4\sin^2 2\theta\cdot  \alpha_1 + \alpha_3-\cos 2\theta\cdot  \alpha_5) (-1-i) D} 
\nonumber \\
&\qquad \overline{+ (-\cos 2\theta\cdot  \alpha_2+\alpha_4-\cos 2\theta\cdot  \alpha_6) (1-i) D }, \nonumber \\
\label{Nayatani-Shoda-eq7} &\left(\alpha_1 + \frac{3}{4} \alpha_3 \right) (-1+i) B
+\left(-\frac{3}{2} \cos 2\theta\cdot  \alpha_2 +\frac{\alpha_4}{4} \right) (1+i) B \\ 
&\qquad +  (-\cos 2\theta\cdot  \alpha_3 +\alpha_5 ) (1-i) D \nonumber \\ 
&\qquad +((-5+6 \cos^2 2\theta)\alpha_2 -\cos 2\theta\cdot \alpha_4 + \alpha_6) (-1-i) D 
\nonumber 
\end{align}
\begin{align}
&= 
\overline{\left(-\cos 2\theta\cdot  \alpha_1 +\frac{\alpha_5}{4} \right) (-1+i) B 
+ \left(\frac{\alpha_2}{4}+\frac{3}{4} \alpha_6 \right) (1+i) B} \nonumber \\
&\qquad \overline{+ (-4\sin^2 2\theta\cdot  \alpha_1 + \alpha_3-\cos 2\theta\cdot  \alpha_5) (1-i) D} 
\nonumber \\
&\qquad \overline{+ (-\cos 2\theta\cdot  \alpha_2+\alpha_4-\cos 2\theta\cdot \alpha_6) (-1-i) D }. \nonumber  
\end{align}

Likewise, one expresses the latter relation of \eqref{period-new} for the homology basis as 
\begin{align}
\label{Nayatani-Shoda-eq2} &\Re\bigg[\left(\frac{3}{4}\alpha_2 + \frac{\alpha_6}{4} \right) A
+ \left(\alpha_1 + \frac{\alpha_3}{4}\right) A \\
&\quad  + (-\cos 2\theta\cdot  \alpha_2+\alpha_4-\cos 2\theta\cdot \alpha_6)C 
+ (-\cos 2\theta\cdot  \alpha_3+\alpha_5) C \bigg]=0,\nonumber \\
\label{Nayatani-Shoda-eq4} &\Re\bigg[\left(\frac{3}{4}\alpha_2 + \frac{\alpha_6}{4} \right) i A
+ \left(\alpha_1 + \frac{\alpha_3}{4}\right) (-i A) \\
&\quad  + (-\cos 2\theta\cdot  \alpha_2+\alpha_4-\cos 2\theta\cdot \alpha_6) i C 
+ (-\cos 2\theta\cdot  \alpha_3+\alpha_5) (-i C )\bigg]=0,\nonumber \\
\label{Nayatani-Shoda-eq6} &\Re\bigg[\left(\frac{3}{4}\alpha_2 + \frac{\alpha_6}{4} \right) (1+i) B
+ \left(\alpha_1 + \frac{\alpha_3}{4}\right) (-1+i) B \\
&\qquad  + (-\cos 2\theta\cdot  \alpha_2+\alpha_4-\cos 2\theta\cdot \alpha_6)(-1-i) D \nonumber \\
&\qquad + (-\cos 2\theta\cdot  \alpha_3+\alpha_5) (1-i) D \bigg]=0,\nonumber \\
\label{Nayatani-Shoda-eq8} &\Re\bigg[\left(\frac{3}{4}\alpha_2 + \frac{\alpha_6}{4} \right) (-1+i) B
+ \left(\alpha_1 + \frac{\alpha_3}{4}\right) (1+i) B \\
&\qquad  + (-\cos 2\theta\cdot  \alpha_2+\alpha_4-\cos 2\theta\cdot \alpha_6)(1-i) D \nonumber \\
&\qquad + (-\cos 2\theta\cdot  \alpha_3+\alpha_5) (-1-i) D\bigg]=0. \nonumber 
\end{align}

\eqref{Nayatani-Shoda-eq1}, \eqref{Nayatani-Shoda-eq3} are equivalent to 
\begin{align}
\label{Nayatani-Shoda-eq9} &\left(\alpha_1 + \frac{3}{4} \alpha_3 \right) A
+  (-\cos 2\theta\cdot  \alpha_3 +\alpha_5 ) C \\
&= 
\overline{
\left(\frac{\alpha_2}{4}+\frac{3}{4} \alpha_6 \right) A 
+ (-\cos 2\theta\cdot  \alpha_2+\alpha_4-\cos 2\theta\cdot \alpha_6) C}, 
\nonumber \\
\label{Nayatani-Shoda-eq10} &\left(-\frac{3}{2} \cos 2\theta\cdot  \alpha_2 +\frac{\alpha_4}{4} \right) A 
+((-5+6 \cos^2 2\theta)\alpha_2 -\cos 2\theta\cdot \alpha_4 + \alpha_6) C \\
&= 
\overline{\left(-\cos 2\theta\cdot  \alpha_1 +\frac{\alpha_5}{4} \right) A 
+ (-4\sin^2 2\theta\cdot  \alpha_1 + \alpha_3-\cos 2\theta\cdot \alpha_5) C }.
\nonumber  
\end{align}
\eqref{Nayatani-Shoda-eq5}, \eqref{Nayatani-Shoda-eq7} are euivalent to 
\begin{align}
\label{Nayatani-Shoda-eq14} &\left(\alpha_1 + \frac{3}{4} \alpha_3 \right) B
+  (-\cos 2\theta\cdot  \alpha_3 +\alpha_5 ) (-D)  
\end{align}
\begin{align}
&= 
-\overline{\bigg\{\left(\frac{\alpha_2}{4}+\frac{3}{4} \alpha_6 \right) B
+ (-\cos 2\theta\cdot  \alpha_2+\alpha_4-\cos 2\theta\cdot \alpha_6) (-D) \bigg\}},\nonumber \\
\label{Nayatani-Shoda-eq15} &\left(-\frac{3}{2} \cos 2\theta\cdot  \alpha_2 +\frac{\alpha_4}{4} \right) B 
+((-5+6 \cos^2 2\theta)\alpha_2 -\cos 2\theta\cdot \alpha_4 + \alpha_6) (-D) \\
&= 
-\overline{\bigg\{\left(-\cos 2\theta\cdot  \alpha_1 +\frac{\alpha_5}{4} \right) B 
+ (-4\sin^2 2\theta\cdot  \alpha_1 + \alpha_3-\cos 2\theta\cdot \alpha_5) (-D)\bigg\}}. \nonumber 
\end{align}
\eqref{Nayatani-Shoda-eq2}, \eqref{Nayatani-Shoda-eq4} are equivalent to 
\begin{align}
\label{Nayatani-Shoda-eq11} & \left(\alpha_1 + \frac{\alpha_3}{4}\right) A 
+ (-\cos 2\theta\cdot  \alpha_3+\alpha_5) C \\
&= 
-\overline{\bigg\{ 
\left(\frac{3}{4}\alpha_2 + \frac{\alpha_6}{4} \right) A
+ (-\cos 2\theta\cdot  \alpha_2+\alpha_4-\cos 2\theta\cdot \alpha_6)C \bigg\}}.
\nonumber  
\end{align}
\eqref{Nayatani-Shoda-eq6}, \eqref{Nayatani-Shoda-eq8} are equivalent to 
\begin{align}
\label{Nayatani-Shoda-eq18} &\left(\alpha_1 + \frac{\alpha_3}{4}\right) B 
+ (-\cos 2\theta\cdot  \alpha_3+\alpha_5) (-D) \\
&= 
\overline{ \left(\frac{3}{4}\alpha_2 + \frac{\alpha_6}{4} \right) B
+ (-\cos 2\theta\cdot  \alpha_2+\alpha_4-\cos 2\theta\cdot \alpha_6) (-D) }.
\nonumber 
\end{align}

The equations \eqref{Nayatani-Shoda-eq9}--\eqref{Nayatani-Shoda-eq18}
are summarized as 
{\small
\begin{align}
\begin{pmatrix}
X_1 & X_2 
\end{pmatrix}
\begin{pmatrix}
\alpha_1 \\
\alpha_5 \\
\overline{\alpha_2} \\
\overline{\alpha_4} \\
\alpha_3 \\
\overline{\alpha_6}
\end{pmatrix}
&=
\begin{pmatrix}
0 \\
0 \\
0 \\
0 \\
0 \\
0
\end{pmatrix}, \label{Nayatani-Shoda-key-relation1}
\end{align}
where 
\begin{align*}
&X_1=\begin{pmatrix}
A & C & \frac{3}{4} A-C \cos 2\theta\\
B & -D & -\frac{3}{4} B-D \cos 2\theta \\
B & -D & \frac{B}{4}+D\cos 2\theta  \\
A & C & -\frac{A}{4}+C \cos 2\theta\\
- \left( A \cos 2\theta +4 C \sin^2 2\theta \right) & \frac{A}{4}-C \cos 2\theta 
& \frac{3}{2} A \cos 2\theta + (5-6 \cos^2 2\theta) C\\
B \cos 2\theta-4 D \sin^2 2\theta & - \left(\frac{B}{4}+D\cos 2\theta\right) 
&\frac{3}{2} B \cos 2\theta + (-5+6 \cos^2 2\theta)D
\end{pmatrix},\\
&X_2=\begin{pmatrix}
C & \frac{A}{4}-C \cos 2\theta & \frac{A}{4}-C \cos 2\theta \\
D & \frac{B}{4}+D \cos 2\theta & -\frac{B}{4}-D \cos 2\theta \\
-D & \frac{3}{4} B +D \cos 2\theta & \frac{3}{4} B+D \cos 2\theta \\
-C & \frac{3}{4} A-C \cos 2\theta & -\frac{3}{4} A + C \cos 2\theta \\
-\frac{A}{4}+C\cos 2\theta & C & -C \\
-\frac{B}{4} -D\cos 2\theta & D & D 
\end{pmatrix}. 
\end{align*}
}

By applying elementary transformations as listed in Appendix A, 
it can be verified that the above system of linear equations is equivalent to 
{\small
\begin{align}
\begin{pmatrix}
Y_1 & Y_2 & Y_3
\end{pmatrix}
\begin{pmatrix}
\alpha_1 \\
\alpha_5 \\
\overline{\alpha_2} \\
\overline{\alpha_4} \\
\alpha_3 \\
\overline{\alpha_6}
\end{pmatrix}
&=
\begin{pmatrix}
0 \\
0 \\
0 \\
0 \\
0 \\
0
\end{pmatrix}, 
\end{align}
where  
\begin{align*}
&Y_1=\begin{pmatrix}
A (AD+BC)^2 & 0 & 0 & 0 \\
0 & -(AD+BC)^2 & 0 & 0 \\
0 & 0 & AD+BC & 0 \\
0 & 0 & 0 & -2C(AD+BC) \\
0 & 0 & 0 & 0 \\
0 & 0 & 0 & 0 
\end{pmatrix},\\
&Y_2=\begin{pmatrix}
\frac{1}{2}A (AD+BC)^2+ \frac{1}{8}A (-AD+BC)^2 \\
(AD+BC)^2 \cos 2\theta+ \frac{1}{4}AB(-AD+BC) \\
\frac{1}{2}(-AD+BC) \\
(AB+(AD-BC)\cos 2\theta) C \\
\frac{1}{16}AC(3AD+BC)[-B (A^2+16 C^2\sin^2 2\theta) +8 (AD+BC)(A \cos 2\theta + 4 C \sin^2 2\theta)] \\
-\frac{1}{16}B D(AD+3 BC)[A (B^2+16 D^2\sin^2 2\theta) +8 (AD+BC)(B \cos 2\theta - 4 D \sin^2 2\theta)]
\end{pmatrix}, \\
&Y_3=\begin{pmatrix}
\frac{1}{2}A(AD+BC) (-AD+BC) \\
AB(AD+BC) \\
AD+BC \\
0 \\
\frac{1}{4}AC(AD+BC)[B (A^2+16 C^2\sin^2 2\theta)-8 (AD+BC)(A \cos 2\theta + 4 C \sin^2 2\theta)] \\
-\frac{1}{4}B D(AD+BC)[A (B^2+16 D^2\sin^2 2\theta) +8 (AD+BC)(B \cos 2\theta - 4 D \sin^2 2\theta)]
\end{pmatrix}. 
\end{align*}
}

It is easy to see that this system has a nontrivial solution if and only if 
the matrix 
$$
\begin{pmatrix} (Y_2)_5 & (Y_3)_5\\ (Y_2)_6 & (Y_3)_6 \end{pmatrix}
$$
is not invertible, where $(Y_i)_j$ is the $j$-th component of $Y_i$. 
In conclusion, the necessary and sufficient condition that 
\eqref{Nayatani-Shoda-key-relation1} has a nontrivial solution is 
that either 
\begin{equation}
A (B^2+16D^2\sin^2 2\theta)+ 8 (AD+BC)(B\cos 2\theta-4D\sin^2 2\theta)=0 
\label{nayatani-shoda-equivalence1}
\end{equation}
or
\begin{equation}
B (A^2+16C^2\sin^2 2\theta)- 8 (AD+BC)(A\cos 2\theta+4C\sin^2 2\theta)=0 
\label{nayatani-shoda-equivalence2}
\end{equation}
holds. 

One can verify that the equation \eqref{nayatani-shoda-equivalence1} 
has a unique solution $\theta_1\approx 0.65 < \pi/4$ 
in the range $0<\theta<\pi/2$. 
We shall give a proof of this fact in Appendix B. 
Note that the change of variable $\theta \mapsto \pi/2-\theta$ transforms 
\eqref{nayatani-shoda-equivalence1} to 
\eqref{nayatani-shoda-equivalence2} and vice versa. 
Therefore, $\theta_2:=\pi/2-\theta_1\approx 0.91 > \pi/4$ gives 
a unique solution of the equation \eqref{nayatani-shoda-equivalence2} 
in the range $0<\theta<\pi/2$. 

If $\theta=\theta_1$, then it is easy to verify that the corresponding 
nontrivial solutions are given by real linear combinations of $\omega_1$ 
and $\omega_2$ as in the statement of Lemma \ref{lemma-extra-ef}.

\section{Proof of Lemma \ref{lemma-symmetry-extra-ef}}

In this section, we shall prove Lemma \ref{lemma-symmetry-extra-ef}. 

Note that $u_i = \langle X_{\omega_i}, N \rangle$, 
where $N$ is the unit normal vector field of $X_{\omega_i}$, 
related to $g_{\theta_1}$ by 
$$
N = {}^t\! \left( 
\frac{2 \Re g_{\theta_1}}{|g_{\theta_1}|^2+1}, 
\frac{2 \Im g_{\theta_1}}{|g_{\theta_1}|^2+1}, 
\frac{|g_{\theta_1}|^2-1}{|g_{\theta_1}|^2+1} 
\right). 
$$

We have $s_1^{*} \omega_1 = \overline{\omega_1}$, 
$s_1^{*} \omega_2 = -\overline{\omega_2}$, 
$$
s_1^{*} \begin{pmatrix} 1-g_{\theta_1}^2 \\ i (1+g_{\theta_1}^2) \\ 2 g_{\theta_1}
\end{pmatrix}
= 
\begin{pmatrix} 1 & 0 & 0 \\ 0 & -1 & 0 \\ 0 & 0 & 1 \end{pmatrix}
\overline{
\begin{pmatrix} 1-g_{\theta_1}^2 \\ i (1+g_{\theta_1}^2) \\ 2 g_{\theta_1} \end{pmatrix}}
,\quad 
s_1^{*} N = \begin{pmatrix} 1 & 0 & 0 \\ 0 & -1 & 0 \\ 0 & 0 & 1 \end{pmatrix} N. 
$$
Since $s_1(p_0)=p_0$, it follows from these formulae that $s_1^{*}u_1=u_1$ 
and $s_1^{*}u_2=-u_2$. 

Let $\psi (z,\,w)=(1/z,\,w/z^3)$. 
By straightforward calculation, we get 
\begin{align*}
& \psi^{*} \left(\frac{dz}{w}\right) =-\frac{z}{w} dz,\quad 
\psi^{*} \left(\frac{dz}{w^3}\right) =-\frac{z^7}{w^3} dz,\quad 
\psi^{*} \left(\frac{z}{w^3}\right) dz =-\frac{z^6}{w^3} dz,\\ 
& \psi^{*} \left(\frac{z^2}{w^3}\right) dz =-\frac{z^5}{w^3} dz,\quad 
\psi^{*} \left(\frac{z^3}{w^3}\right) dz =-\frac{z^4}{w^3} dz,\quad 
\psi^{*} \left(\frac{z^4}{w^3}\right) dz =-\frac{z^3}{w^3} dz. 
\end{align*}
Therefore, 
\begin{align*}
\psi^{*} \omega_1 
&=
z^2\biggl(
\frac{AD+3 BC}{4(AD+BC)} \frac{z^4+2\cos 2\theta\cdot z^2+1}{w^3}dz 
+\frac{AD+3 BC}{4(AD+BC)} \frac{z^5}{w^3} dz - \frac{z^4}{w^3} dz \\
&-\frac{AB+(AD-BC)\cos 2\theta}{2(AD+BC)} \frac{z^3}{w^3} dz 
-\frac{AB+2(AD+BC)\cos 2\theta}{2(AD+BC)} \frac{z^2}{w^3} dz\\ 
&-\frac{3AD+ BC}{4(AD+BC)} \frac{z}{w^3}dz \biggr)\\
&=
z^2\biggl(
\frac{-3AD-BC}{4(AD+BC)} \frac{z^4}{w^3} dz
+ \frac{-AB+(-AD+BC)\cos 2\theta}{2(AD+BC)} \frac{z^2}{w^3}dz\\
& + \frac{AD+3 BC}{4(AD+BC)} \frac{dz}{w^3} 
- \frac{AB+2(AD+BC)\cos 2\theta}{2(AD+BC)} \frac{z^3}{w^3}dz
- \frac{z}{w^3}dz
\end{align*}
\begin{align*}
&+\frac{AD+3BC}{4(AD+BC)} \frac{z^5 + 2\cos 2\theta\cdot  z^3 + z}{w^3}dz \biggr)\\ 
&=-z^2 \omega_1. 
\end{align*}
Likewise, we obtain $\psi^{*} \omega_2=z^2 \omega_2$. 
Since we also have 
$$
\psi^{*}
\begin{pmatrix}
1-g_{\theta_1}^2 \\
i (1+g_{\theta_1}^2) \\
2 g_{\theta_1}
\end{pmatrix}
= 
\frac{1}{z^2} 
\begin{pmatrix}
-1 & 0 & 0 \\
0 & 1 & 0 \\
0 & 0 & 1
\end{pmatrix}
\begin{pmatrix}
1-g_{\theta_1}^2 \\
i (1+g_{\theta_1}^2) \\
2 g_{\theta_1}
\end{pmatrix}, 
\quad \psi^{*} N=
\begin{pmatrix}
1 & 0 & 0 \\
0 & -1 & 0 \\
0 & 0 & -1
\end{pmatrix}
N, 
$$
and $\psi(p_0)=p_0$, we find $\psi^{*} u_1 = u_1$ and $\psi^{*} u_2 = -u_2$. 
Since $s_3=\psi\circ s_1$, we conclude that $s_3^{*} u_1 = u_1$ and 
$s_3^{*} u_2 = u_2$. 

We have $j^{*} \omega_1 = -\omega_1$, $j^{*} \omega_2 = -\omega_2$, 
$$
j^{*} \begin{pmatrix} 1-g_{\theta_1}^2 \\ i (1+g_{\theta_1}^2) \\ 2 g_{\theta_1}
\end{pmatrix}
= 
\begin{pmatrix} 1-g_{\theta_1}^2 \\ i (1+g_{\theta_1}^2) \\ 2 g_{\theta_1} \end{pmatrix}
,\quad 
j^{*} N = N, 
$$
from which it follows that 
\begin{align*}
j^* u_1(p) 
& = \left\langle \Re\int^{j(p)}_{j(p_0)} \begin{pmatrix} 1-g_{\theta_1}^2 \\ 
i (1+g_{\theta_1}^2) \\ 2 g_{\theta_1} \end{pmatrix}\omega_1, N(j(p)) 
\right\rangle 
+ 
\left\langle \Re\int^{j(p_0)}_{p_0} \begin{pmatrix} 1-g_{\theta_1}^2 \\ 
i (1+g_{\theta_1}^2) \\ 2 g_{\theta_1} \end{pmatrix}\omega_1, N(j(p)) 
\right\rangle\\ 
& = - u_1(p) + \langle c_1, N(p) \rangle,
\end{align*}
where $c_1 = \Re\int^{j(p_0)}_{p_0} {}^t (1-g_{\theta_1}^2,  
i (1+g_{\theta_1}^2), 2 g_{\theta_1})\, \omega_1$, and 
$j^*u_2 = -u_2 + \langle c_2, N \rangle$. 

\section*{Appendix A}

As mentioned in the proof of Lemma \ref{lemma-extra-ef}, the system 
\eqref{Nayatani-Shoda-key-relation1} of linear equations 
can be reduced to an equivalent one of simpler form 
by applying elementary transformations. 
For the reader's convenience, we shall list all the 
elementary transformations explicitly. 

We apply the following operations, where ${\rm R}j$ ($1\leq j\leq 6$) denotes the $j$-th row, 
to the matrix $\begin{pmatrix}X_1 & X_2\end{pmatrix}$. 
\begin{enumerate}
\renewcommand{\theenumi}{\roman{enumi}}
\renewcommand{\labelenumi}{(\theenumi)}
\item ${\rm R}4 \longrightarrow {\rm R}4-{\rm R}1$. 
\item ${\rm R}5 \longrightarrow {\rm R}5+{\rm R}1\times \cos 2\theta$. 
\item ${\rm R}3 \longrightarrow {\rm R}3-{\rm R}2$. 
\item ${\rm R}6 \longrightarrow {\rm R}6-{\rm R}2\times \cos 2\theta$. 
\item ${\rm R}5 \longrightarrow {\rm R}5\times A+{\rm R}1\times 4C\sin^2 2\theta$. 
\item ${\rm R}6 \longrightarrow {\rm R}6\times B+{\rm R}2\times 4D\sin^2 2\theta$. 
\item ${\rm R}2 \longrightarrow {\rm R}2+{\rm R}3\times 1/2$. 
\item ${\rm R}6 \longrightarrow {\rm R}6+{\rm R}3\times (-B\cos 2\theta+2D\sin^2 2\theta)$. 
\item ${\rm R}1 \longrightarrow {\rm R}1+{\rm R}4\times 1/2$. 
\item ${\rm R}5 \longrightarrow {\rm R}5+{\rm R}4\times (A\cos 2\theta+2C\sin^2 2\theta)$. 
\item ${\rm R}2 \longrightarrow {\rm R}2\times A+{\rm R}1\times (-B)$. 
\item ${\rm R}5 \longrightarrow {\rm R}5\times C+{\rm R}4\times (-A^2/8)$. 
\item ${\rm R}6 \longrightarrow {\rm R}6\times D+{\rm R}3\times (-B^2/8)$. 
\item ${\rm R}1 \longrightarrow {\rm R}1\times (AD+BC)+{\rm R}2\times C$. 
\item ${\rm R}3 \longrightarrow {\rm R}3\times C+{\rm R}4\times (-D)$. 
\item ${\rm R}5 \longrightarrow {\rm R}5\times (AD+BC)+{\rm R}2\times (A^2 C/4+4C^3\sin^2 2\theta)$. 
\item ${\rm R}6 \longrightarrow {\rm R}6\times (AD+BC)+{\rm R}2\times (-B^2 D/4-4D^3\sin^2 2\theta)$. 
\item ${\rm R}1 \longrightarrow {\rm R}1\times (AD+BC)+{\rm R}3\times A(-AD+BC)/4$. 
\item ${\rm R}2 \longrightarrow {\rm R}2\times (AD+BC)+{\rm R}3\times AB/2$. 
\item ${\rm R}4 \longrightarrow {\rm R}4\times (AD+BC)+{\rm R}3\times (A - 2 C \cos 2\theta)$. 
\item ${\rm R}5 \longrightarrow {\rm R}5\times (AD+BC)+{\rm R}3\times [-A^2 (A^2 D/8+ (AD+B C) C \cos 2\theta + 6 C^2 D \sin^2 2\theta) 
-4 A B C^3\sin^2 2\theta]$. 
\item ${\rm R}6 \longrightarrow {\rm R}6\times (AD+BC)+{\rm R}3\times [-B^2(-B^2 C/8+(A D+BC) D \cos 2\theta - 6 C D^2 \sin^2 2\theta) 
+ 4 A B D^3 \sin^2 2\theta]$. 
\end{enumerate}
Then we finally obtain the matrix 
$\begin{pmatrix} Y_1 & Y_2 & Y_3 \end{pmatrix}$ 
as in the proof of Lemma \ref{lemma-extra-ef}. 

\section*{Appendix B} 
In this appendix, we prove that the equation \eqref{nayatani-shoda-equivalence1} 
has a unique solution $\theta_1 < \pi/4$ in the range $0<\theta<\pi/2$. 

We first prove that \eqref{nayatani-shoda-equivalence1} has a unique 
solution in the range $0<\theta<\pi/4$. 
Though it is possible to verify this fact by a direct elementary argument, 
here we present an indirect one, assuming that \eqref{nayatani-shoda-equivalence1} 
has no solutions in the range $\pi/4\leq \theta<\pi/2$, which we will prove 
afterwards. 
Since the left-hand side of \eqref{nayatani-shoda-equivalence1} is positive 
near $\theta = 0$ and negative at $\theta = \pi/4$, 
\eqref{nayatani-shoda-equivalence1} has at least one solution by the intermediate 
value theorem. 
On the other hand, \eqref{nayatani-shoda-equivalence2} has no solutions in the 
range $0<\theta< \pi/4$ by the remark at the end of the proof of Lemma 
\ref{lemma-extra-ef}. 
Suppose that there is more than one solution of \eqref{nayatani-shoda-equivalence1}, 
and let $\varphi_1<\varphi_2$ be the first and second smallest ones. 
Then, since the argument for proving Theorem \ref{theorem-index} depends only on the fact 
that $\theta_1$ is a solution of \eqref{nayatani-shoda-equivalence1}, we deduce 
that the number of eigenvalues of $-\Delta_\theta$ less than $2$ decreases by two 
each time when $\theta$ passes $\varphi_1$ and $\varphi_2$. 
But this is impossible because there are exactly three such eigenvalues for 
$\theta< \varphi_1$. 
Thus, the solutions of \eqref{nayatani-shoda-equivalence1} must be unique. 

We now proceed to prove that the equation \eqref{nayatani-shoda-equivalence1} 
has no solutions in the range $\pi/4\leq\theta<\pi/2$. 
We start by rewriting the integrals $A, B, C, D$ using the complete elliptic integrals 
$$
K(\textbf{k})=\int_{0}^{\frac{\pi}{2}} 
\frac{d\theta}{\sqrt{1-\textbf{k}^2\sin^2 \theta}},\quad
E(\textbf{k})=\int_{0}^{\frac{\pi}{2}} \sqrt{1-\textbf{k}^2\sin^2 \theta}
\,d\theta,  
$$
defined for $0<\textbf{k}<1$. 
Clearly, $K(\textbf{k})$ (resp. $E(\textbf{k})$) is a monotone increasing 
(resp. decreasing) function of $\textbf{k}$. 
Computing with the change of variable $u=\sqrt{t}-1/\sqrt{t}$ and using \textbf{222} 
of \cite{ByrdFriedman}, we obtain 
\begin{align*}
&A = \frac{2}{\sqrt{2(1+\sin \theta)}} K(k),\quad 
B = \frac{2}{\sqrt{2(1+\cos \theta)}} K(l), \\
&C = \frac{1}{4\sqrt{2(1+\sin \theta)} \sin^2\theta (1-\sin\theta)} 
\left(E(k)-(1-\sin\theta) K(k)\right),\\ 
&D = \frac{1}{4\sqrt{2(1+\cos \theta)} \cos^2\theta (1-\cos\theta)} 
\left(E(l)-(1-\cos\theta) K(l)\right),  
\end{align*}
where $k = \sqrt{2\sin\theta/(1+\sin\theta)}$ and $l = \sqrt{2\cos\theta/(1+\cos\theta)}$. 

The left-hand side of \eqref{nayatani-shoda-equivalence1} can be rewritten as 
$$
\cos 2\theta (AB^2 \cos 2\theta +8 A B D +8 B^2 C)+\sin^2 2\theta (AB^2-16AD^2-32BCD). 
$$
Therefore, it suffices to verify that both 
\begin{equation}\label{two-inequalities}
AB^2 \cos 2\theta +8 A B D +8 B^2 C>0,\quad AB^2-16AD^2-32BCD<0
\end{equation}
hold in the range $\pi/4\leq \theta <\pi/2$. 

We first reduce these inequalities to several simpler ones, with details 
discussed later on. 
The former inequality of \eqref{two-inequalities} follows from 
\begin{equation}\label{NS-estimate1}
A\cos 2\theta +8 C > 0,\quad \pi/4\leq \theta <\pi/2. 
\end{equation}
For the latter inequality of \eqref{two-inequalities}, since 
\begin{eqnarray*}
\lefteqn{A B^2-16A D^2-32 B C D}\\ 
&=& \left\{\begin{array}{l}
A\left( B-\frac{192}{25} D \right) \left( B+\frac{25}{12} D \right)
+ BD \left( \frac{1679}{300}A-32C \right),\quad \pi/4\leq \theta\leq 5\pi/16,\\ 
A (B - 10 D) \left( B + \frac{8}{5} D \right)
+ BD \left( \frac{42}{5}A - 32 C \right),\quad 5\pi/16\leq \theta\leq 3\pi/8,\\ 
A(B-16D)(B+D) + BD(15A-32C),\quad 3\pi/8\leq \theta<\pi/2,
\end{array}\right. 
\end{eqnarray*}
it suffices to show 
\begin{align}
\label{NS-estimate2} & 25 B-192 D<0,\quad \pi/4\leq \theta \leq 5\pi/16,\\ 
\label{NS-estimate3} & 1679 A-9600 C < 0,\quad \pi/4\leq \theta \leq 5\pi/16, \\ 
\label{NS-estimate4} & B-10D<0,\quad 5\pi/16 \leq \theta \leq 3\pi/8, \\ 
\label{NS-estimate5} & 21A - 80 C<0,\quad 5\pi/16 \leq \theta \leq 3\pi/8,\\  
\label{NS-estimate6} & B-16D<0,\quad 3\pi/8\leq \theta <\pi/2,\\ 
\label{NS-estimate7} & 15A-32C<0,\quad 3\pi/8\leq \theta <\pi/2.
\end{align}

We now present a detailed proof of \eqref{NS-estimate2}. 
Since the proofs of \eqref{NS-estimate1} and 
\eqref{NS-estimate3}--\eqref{NS-estimate7} are similar, 
they are left to the reader.  
We have  
$$
25 B - 192 D = f(l)\, [ (49l^4-96l^2+96)\, (1-l^2)\, K(l) - 12\, (2-l^2)^3\, E(l) ], 
$$
where $f(l)$ is a positive function of $l$.  
Therefore, one must show that 
$$
(49l^4-96l^2+96)(1-l^2) K(l)-12(2-l^2)^3E(l)<0
$$
in the range 
$$
0.7142\cdots = \frac{2\cos\frac{5}{16} \pi}{1+\cos \frac{5}{16}\pi}
\leq l^2\leq
\frac{2\cos\frac{\pi}{4}}{1+\cos \frac{\pi}{4}} = 0.8284\cdots. 
$$
Using 
$$
\frac{d}{d\textbf{k}} [ (1-\textbf{k}^2) K(\textbf{k}) ]
= \frac{E(\textbf{k})}{\textbf{k}}-\frac{1+\textbf{k}^2}{\textbf{k}} 
K(\textbf{k}),\quad 
\frac{d}{d\textbf{k}} E(\textbf{k})= \frac{E(\textbf{k})-K(\textbf{k})} 
{\textbf{k}}
$$
(cf.~ \cite[\textbf{710}]{ByrdFriedman}), 
we obtain 
\begin{eqnarray}\label{derivative11}
\lefteqn{\frac{d}{dl} [ (49l^4-96l^2+96)(1-l^2) K(l)-12(2-l^2)^3E(l) ]}\\ 
&=& l [ -(257 l^4-507 l^2+336 )K(l)+(84 l^4-311 l^2+336) E(l) ]. \nonumber 
\end{eqnarray}
Observe that $257 l^4-507 l^2+336$ and $84 l^4-311 l^2+336$ are positive and 
monotone decreasing in the range $0.71< l^2< 0.83$. 
Then we can show that the right-hand side of \eqref{derivative11} is negative 
in the range $0.71< l^2< 0.83$ by estimating it in $0.71< l^2\leq 0.81$ and $0.81\leq l^2< 0.83$ 
separately. 
E.g., in $0.71< l^2\leq 0.81$, 
\begin{eqnarray*}
\lefteqn{-(257 l^4-507 l^2+336 )K(l)+(84 l^4-311 l^2+336) E(l)} \\
&\leq& -(257\cdot 0.81^2-507 \cdot 0.81+336 )K\left(\sqrt{0.71}\right) \\
&&+(84 \cdot 0.71^2-311 \cdot 0.71+336) E\left(\sqrt{0.71}\right) \\
&=& -1.723\cdots < 0. 
\end{eqnarray*} 
Therefore, $(49l^4-96l^2+96) (1-l^2) K(l) - 12 (2-l^2)^3 E(l)$
is monotone decreasing. 
Since its value at $l^2=0.714$ is $-0.033\cdots
<0$, \eqref{NS-estimate2} is proved.

\end{document}